\documentclass[12pt]{article}
\usepackage{graphicx} 
\usepackage{epstopdf}
\usepackage{subfigure}
\usepackage{color}
\usepackage[english]{babel}
\usepackage{amsmath,amsthm}
\usepackage{amssymb}
\usepackage{amsfonts}
\usepackage{booktabs}

\usepackage{graphicx}
\usepackage{subfigure}
\usepackage{tabularx}
\title{On the dispersability of graph bundles over cycles}

\author { Zeling Shao, Xiaoxiang Yu, Zhiguo Li{$^*$}\\
{\small School of Science, Hebei University of Technology, Tianjin 300401, China}
\date{}
\footnote{Corresponding author. E-mail: zhiguolee@hebut.edu.cn}
}

\begin{document}
\baselineskip 0.65cm

\maketitle

\begin{abstract}

  In this paper, the dispersability of the Cartesian graph bundle over two cycles is completely solved. We show the Cartesian graph bundle $G$ over two cycles is dispersable if $G$ is bipartite; otherwise, $G$ is nearly dispersable.



\noindent\textbf{2020 MR Subject Classification.} 05C10
\end{abstract}

\section{Introduction}
The research of book embedding is of great importance since it has applications in several areas of computer science, such as fault-tolerant computing, sorting with parallel stacks, and so on(see $[1,2,3]$). Let $\psi$ be a permutation of all vertices of a graph $G$.~A \emph{layout} $\Psi=(G,\psi)$ for $G$ is to arrange all vertices along a circle in the order $\psi$ and join the edges of $G$ as chords.~Let $S$ be a color set and $|S|=m$.~A triple $(G,\psi,c)$ is an $m$-page \emph{book embedding} if $c:E(G) \rightarrow S$ is an edge-coloring such that $c(e^{\prime})\neq c(e^{\prime\prime})$ when $e^{\prime}$ and $e^{\prime\prime}$ cross in $\Psi$.~The \emph{book thickness}~$bt(G)$ of $G$ is the minimum integer $m$ such that an $m$-page book embedding exists.~A book embedding $(G,\psi,c)$ is \emph{matching} if the edge-coloring $c$ is proper.~The \emph{matching book thickness} $mbt(G)$ of $G$ is the minimum integer $m$ such that an $m$-page matching book embedding exists.~We call $G$ \emph{dispersable} if $mbt(G)=\Delta(G)$ and \emph{nearly dispersable} if $mbt(G)=\Delta(G)+1$.\\
\indent The dispersability of some families of graphs has been studied. The dispersability of complete bipartite graphs, even cycles, binary $n$-cubes$(n \geq 1)$, trees, and complete graphs have been solved in $[3],[5]$. Shao et al. have obtained the dispersability of the generalized Petersen graph $[6]$ and the pseudo-Halin graph $[7]$. The dispersability of the Cartesian product of two cycles has been solved. Kainen ${[8]}$ showed that $C_{2n}\Box C_{2m}$ is dispersable and  $C_{2n}\Box C_{2m+1}$ is nearly dispersable. Shao, Liu, Li $[9]$ showed that $mbt(K_{n} \Box C_{q})$$=$$\Delta (K_{n}\Box C_{q})+1$ for $n,q \geq 3$,~which implies $C_{3} \Box C_{q}$ is nearly dispersable.~In $[4]$ it was proved that $C_{5}\Box C_{n}$ is nearly dispersable for $n\geq 5$.~The authors $[10]$ showed $C_{2m+1}\Box C_{2n+1}$ is nearly dispersable, $m,n \geq 3$, which completely solved the dispersability of the Cartesian product of two cycles.\\
\indent Circulant graphs have deserved significant attention in the last decades either theoretically or through their applications on building interconnection networks for parallel computing.~Let us review some definitions and notations.~For convenience, unless otherwise specified,~let the group $\mathbb{Z}_{n}$ be the set $\{1,2,\cdots,n\}$ throughout this paper.~Let $S$ be a subset of $\mathbb{Z}_{n}$ such that any $k\in S$ satisfies $1\leq k\leq \lfloor n/2 \rfloor$.~A \emph{circulant graph} $C(\mathbb{Z}_{n},S)$ is the graph whose vertex set is $\mathbb{Z}_{n}$,~two vertices $i,j$ are adjacent if $i-j\equiv k\pmod{n}(k\in S)$.~Each element of $S$ is called a jump length.\\
\indent Joslin, Kainen and Overbay $[11]$ gave some results for the dispersability of the circulant graph $C(\mathbb{Z}_{n},S)$, where $S$ is a subset of $\{1,2,3\}$ with size $2$. In addition, they showed if $n$ is a multiple of $2k+1$, then the circulant graph $C(\mathbb{Z}_{n},\{1,2,\cdots,k\})$ is nearly dispersable; if $n$ is a multiple of $12$, then the circulant graph $C(\mathbb{Z}_{n},\{1,2,3\})$ is nearly dispersable; if $2k~|~n~(k\geq 3)$, then $C(\mathbb{Z}_{n},\{1,3,5,\cdots,2d+1\})$ is dispersable, where $2d+1$ is the largest odd integer not exceeding $k$.

In [12],  the  classification  of circulant graphs $C(\mathbb{Z}_{n},\{k_1,k_2\})$ is obtained and the dispersability of circulant graphs $C(\mathbb{Z}_{n},\{k_1,k_2\})$ is all solved except the case with the Cartesian graph bundle over two cycles as its component.



Graph bundle is a generalization of covering graphs and Cartesian products of graphs. The notion follows the definition of fiber bundles and vector
bundles in topology as a space which locally looks like a product. Let us review some definitions and properties of graph bundle $[13]$. Let $B,F$ be graphs,~a graph $G$ is a \emph{Cartesian graph bundle} with the fiber $F$ over the base graph $B$ if there is a graph map $p: G \rightarrow B$ such that for each vertex $v\in V(B)$, $p^{-1}(v)\cong F,$ and for each edge $e\in E(B)$, $p^{-1}(e)\cong K_2 \Box F.$ Let $\varphi: E(B)\rightarrow Aut(F)$ be a mapping which assigns an automorphism of the graph $F$ to any edge of $B$. The bundle $G$ is denoted by $G=B\Box^{\varphi}F$. 

An automorphism $\varphi$ of a cycle $C_t$ is of two types [14]. A cyclic shift of the cycle $C_t$ by $d$ elements is called the cyclic $d$-shift and another type of $C_t$ is called a reflection, $d\in \mathbb{Z}_{t}$. As for $\varphi$ is a reflection, $\varphi$ has one fixed point if $t$ is odd; otherwise, there are two isomorphism classes of $\varphi$, depending on the number of fixed points is 0 or 2. For convenience, we assume that the vertex set of $C_s \Box^{\varphi} C_t$ is $\{(p,q)~|~p\in \mathbb{Z}_{s},q\in \mathbb{Z}_{t}\}$ (see Fig.1 for the graph bundle $C_5 \Box^{\varphi} C_6$, where $\varphi$ (left) is a 1-shift, and $\varphi$ (right) is a reflection without fixed points).

\vspace{-0.6em}
\begin{figure}[htbp]
\centering
\includegraphics[height=3.6cm, width=0.66\textwidth]{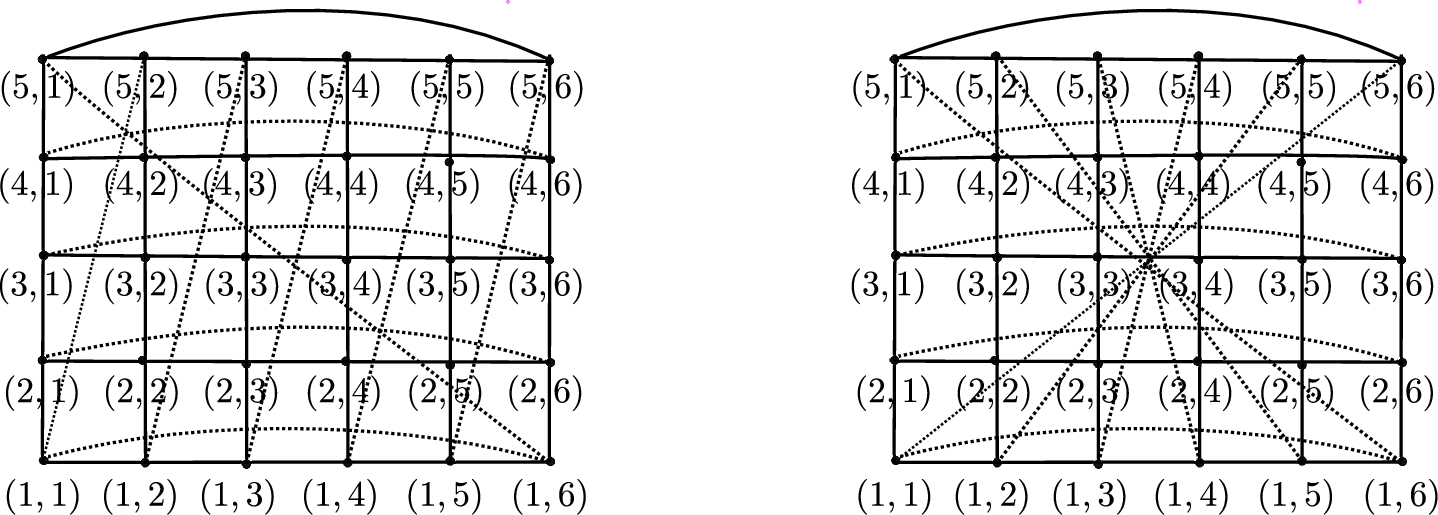}\\
\small Fig.1~~The graph bundle $C_5 \Box^{\varphi} C_6$, where $\varphi$ (left) is a 1-shift, and $\varphi$ (right) is a reflection without fixed points.
\end{figure}
\vspace{-1.2em}
\indent In this paper, we consider the dispersability of the graph bundle $G=C_s \Box^{\varphi} C_t$ and get that $G$ is dispersable if it is bipartite; otherwise, $G$ is nearly dispersable.\\
\indent The paper is organized as follows.~We introduce some properties of matching book embedding and the graph bundle $G=C_s \Box^{\varphi} C_t$ in section 2. For the case where $\varphi$ is a nontrivial shift, section 3 gives the dispersability of the graph bundle $G$. Section 4 discusses the case where $\varphi$ is a reflection. Finally, we make a conclusion in Section 5.

\section{Preliminaries}
Let $gcd(n,k)$ be the greatest common denominator of $n$ and $k$. Denote the reverse and the size of an ordered set $X$ by $X^{-}$ and $|X|$, respectively. 
Firstly, we recall some preliminary properties of matching book embedding.

\noindent \textbf{Definition 2.1.~}\emph{Given a layout $\Psi=(G,\omega)$ of a graph $G$, let $S$ be a color set of size $mbt(G)$.}
\emph{We say the graph $G$ can be colored well in the order $\omega$ if there is a proper edge coloring $c:E(G)\rightarrow S$ such that $c(e)\neq c(e^{\prime})$ when $e$ and $e^{\prime}$ cross in $\Psi$ for $e,e^{\prime}\in E(G)$.}

\noindent \textbf{Note.~}If $G$ is colored well, obviously, an $mbt(G)$-page matching book embedding of $G$ is obtained.


\vspace{0.0em}
\noindent \textbf{Lemma 2.1.~}$^{[5]}$\emph{~If a regular graph $G$ is dispersable, then $G$ is bipartite.}

\vspace{0.0em}
\noindent \textbf{Lemma 2.2.~}$^{[5]}$\emph{~For any simple graph $G$, we have $\Delta(G)\leq \chi^{\prime}(G) \leq mbt(G)$, where $\chi^{\prime}(G)$ is the edge chromatic number of $G$.}

\vspace{0.2em}
Next, we list some properties of the Cartesian graph bundle $C_s \Box^{\varphi} C_t$, where $s,t\geq 3$ and $\varphi$ is an automorphism of $C_t$.

\vspace{0.0em}
\noindent \textbf{Fact 2.1.~}\emph{The Cartesian graph bundle $C_s \Box^{\varphi} C_t$ is $4$-regular.}

A nonempty graph is 2-chromatic if and only if it is bipartite. It is easy to check the following result holds by the chromatic number of $C_s \Box^{\varphi} C_t$ in $[14]$.

\vspace{0.0em}
\noindent \textbf{Lemma 2.3.~}\emph{$(i)$~Let $d \in \mathbb{Z}_{t}$. If $\varphi$ is a nontrivial $d$-shift, then $C_s \Box^{\varphi} C_t$ is bipartite if and only if $t$ is even and $s,d$ have the same parity; otherwise, $C_s \Box^{\varphi} C_t$ is nonbipartite.}

\noindent\emph{$(ii)$~If $\varphi$ is a reflection, then $C_s \Box^{\varphi} C_t$ is bipartite if and only if $\varphi$ has no fixed points, $s$ is odd or $\varphi$ has two fixed points, $s$ is even; otherwise, $C_s \Box^{\varphi} C_t$ is nonbipartite.}

\vspace{0.0em}
\noindent \textbf{Lemma 2.4.~}$^{[14]}$\emph{~If $d_1+d_2\equiv 0\pmod{t}$, then $C_s \Box^{\varphi_1} C_t\cong C_s \Box^{\varphi_2} C_t$, where $\varphi_1$ and $\varphi_2$ are the $d_1$-shift and $d_2$-shift of the cycle $C_t$, respectively.}

\vspace{0.2em}
For future reference, we state as some lemmas the observation that the Cartesian graph bundle $C_s \Box^{\varphi} C_t$ admits a cycle decomposition.  In addition, we introduce some notations, which plays an important role to give vertex orderings and edge-colorings of the optimal matching book embedding of $C_s \Box^{\varphi} C_t$ in Section 3 and Section 4.\\
\indent We assume that the vertex set of the base $C_s$ is $\mathbb{Z}_{s}$. Let $F_i=p^{-}(v_i)$ be the $t$-cycle induced by the fiber $C_t$, where $V({F_i})=\{(i,1),(i,2),\cdots,(i,t)\}$ and vertices $(i,m),(i,n)$ are adjacent if $m-n\equiv 1\pmod {t}$, $i\in V(C_{s})$. Denote the clockwise cyclic ordered vertex set of the cycle ${F_i}$ by $A_i$, $i\in\mathbb{Z}_{s}$. Specifically, $A_i=\{(i,1),(i,2),\cdots,(i,t)\}$. Similarly, let $B_j^{\prime}=\{(1,j),(2,j),\cdots,(s,j)\}$ be a vertex set, and $B_j$ be a permutation of $B_j^{\prime}$, where $B_j^{\prime}=\{(1,j),(2,j),\cdots,(s,j)\}$, $j\in\mathbb{Z}_{t}$.\\
\noindent \textbf{Lemma 2.5.~}$^{[12]}~$\emph{If $d\neq n/2$, then $C(\mathbb{Z}_{t},\{d\})$ is the edge disjoint union of $gcd(t,d)$ cycles of length $t/gcd(t,d)$.}\\
\indent Let $X$ be a proper subset of $V(G)$. To \emph{shrink} $X$ is to delete all edges between vertices of $X$ and then identify the vertices of $X$ into a single vertex.\\
\noindent \textbf{Lemma 2.6.~}\emph{Let $\varphi$ be a nontrivial $d$-shift of the fiber $C_t$, $d \in \mathbb{Z}_{t}$. If $G^{\prime}$ is the spanning subgraph of $C_s \Box^{\varphi} C_t$ by deleting all of these $t$-cycles ${F_1,F_2,\cdots,F_s}$, then $G^{\prime}$ has $gcd(t,d)$ isomorphic connected components $H_1,H_2,\cdots,H_{gcd(t,d)}$ and each one of these components is a cycle of length $st/gcd(t,d)$.}\\
\noindent \textbf{Proof.~}If $d=t/2$, we have $V(H_k)=B^{\prime}_{k}\cup B^{\prime}_{k+d}$, $k\in\mathbb{Z}_{gcd(t,d)}$ (see Fig.2(left) for $C_5 \Box^{\varphi_{1}} C_8$, where $\varphi_{1}$ is a 4-shift). If $d\neq t/2$,~let $H^{\prime}$ be a graph which is obtained from $G^{\prime}$ by shrinking $B_{j}^{\prime}$ into a single vertex $j$ for any $j\in\mathbb{Z}_{t}$.~It is easy to see $H^{\prime}$ is a circulant graph $C(\mathbb{Z}_{t},\{d\})$. By Lemma 2.5, $C(\mathbb{Z}_{t},\{d\})\cong H_1^{\prime}\cup H_2^{\prime}\cup\cdots\cup H_{gcd(t,d)}^{\prime}$, where $H_k^{\prime}$ is a cycle of length $t/gcd(t,d)$, $k\in\mathbb{Z}_{gcd(t,d)}$. Specifically, $V(H_k^{\prime})=\{k+(l-1)d~|~l\in\{1,2,\cdots,t/gcd(t,d)\}\}(\text{mod}~t)$. By the inverse operation of shrinking $B_j^{\prime}$, $G^{\prime}$ is the edge disjoint union of the cycles $H_1,H_2,\cdots,H_{gcd(t,d)}$, where $V(H_k)=\underset{l=1}{\overset{t/gcd(t,d)}{\cup}}B^{\prime}_{k+(l-1)d}$ (see Fig.2(right) for $C_5 \Box^{\varphi_{2}} C_7$, where $\varphi_{2}$ is a 3-shift).
\hfill{$\square$}

\vspace{-0.8em}
\begin{figure}[htbp]
\centering
\includegraphics[height=4.3cm, width=0.84\textwidth]{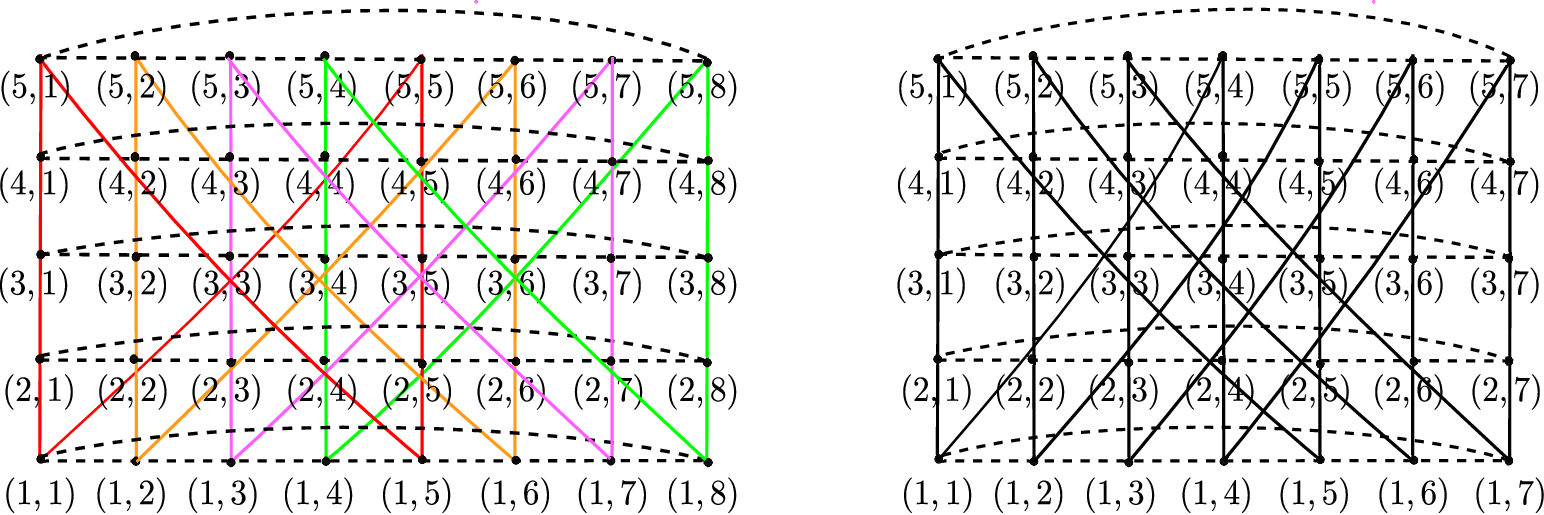}\\
\small Fig.2~~The graph bundle $C_5 \Box^{\varphi_{1}} C_8$, where $\varphi_{1}$ is a 4-shift(left) and $C_5 \Box^{\varphi_{2}} C_7$, where $\varphi_{2}$ is a 3-shift(right).
\end{figure}
\vspace{-0.6em}

Let $V_i$ be the clockwise cyclic ordered vertex set of the cycle ${H_i}$, $1\leq i\leq gcd(t,d)$. Specifically, $V_i=\{(1,i),(2,i),\cdots,(s,i+(t/gcd(t,d)-1)d)\}$.

\vspace{0.1em}
\noindent \textbf{Lemma 2.7.~}\emph{Let $\varphi$ be a reflection of the fiber $C_t$. If $G^{\prime\prime}$ is the spanning subgraph of $C_s \Box^{\varphi} C_t$ by deleting all edges of these $s$ cycles ${F_i}$, $1\leq i\leq s$, we have}\\
\noindent\emph{$(i)$~If $\varphi$ is a reflection without fixed points, then $t$ is even and $G^{\prime\prime}$ is the edge disjoint union of $t/2$ cycles $D_j$, $j\in \mathbb{Z}_{t/2}$.}\\
\noindent\emph{$(ii)$~If $\varphi$ is a reflection with two fixed points, then $t$ is even and $G^{\prime\prime}$ is the edge disjoint union of $t/2+1$ cycles $D_j,j\in \mathbb{Z}_{t/2+1}$.}\\
\noindent\emph{$(iii)$~If $\varphi$ is a reflection with one fixed point, then $t$ is odd and $G^{\prime\prime}$ is the edge disjoint union of $(t+1)/2$ cycles $D_j,j\in \mathbb{Z}_{(t+1)/2}$.}\\
\noindent \textbf{Proof.~}$(i)$~It is easy to see that  $t$ is even and $V(D_j)=B_j^{\prime}\cup B_{t+1-j}^{\prime},j\in \mathbb{Z}_{t/2}$ (see Fig.3a for $C_5 \Box^{\varphi_{1}} C_8$, where $\varphi_{1}$ is a reflection without fixed points).\\
\noindent $(ii)$~Since $\varphi$ has two fixed points, the integer $t$ is even. Without loss of generality, we assume two fixed points of $C_t$ are 1 and $t/2+1$. Thus $V(D_1)=B_1^{\prime}$, $V(D_j)=B_j^{\prime}\cup B_{t+2-j}^{\prime}$, $V(D_{t/2+1})=B^{\prime}_{t/2+1}$, $2\leq j \leq t/2$ (see Fig.3b for $C_5 \Box^{\varphi_{2}} C_8$, where $\varphi_{2}$ is a reflection with two fixed points).\\
\noindent $(iii)$~Obviously, the integer $t$ is odd. Without loss of generality, we assume the fixed point of $C_t$ is $(t+1)/2$. So $V(D_j)=B_j^{\prime}\cup B_{t+1-j}^{\prime}$, $V(D_{(t+1)/2})=B_{(t+1)/2}^{\prime}, j\in \mathbb{Z}_{(t-1)/2}$ (see Fig.3c for $C_5 \Box^{\varphi_{3}} C_9$, where $\varphi_{3}$ is a reflection with one fixed point).
\hfill{$\square$}

\vspace{-0.6em}
\begin{figure}[htbp]
\centering
\includegraphics[height=9.cm, width=0.8\textwidth]{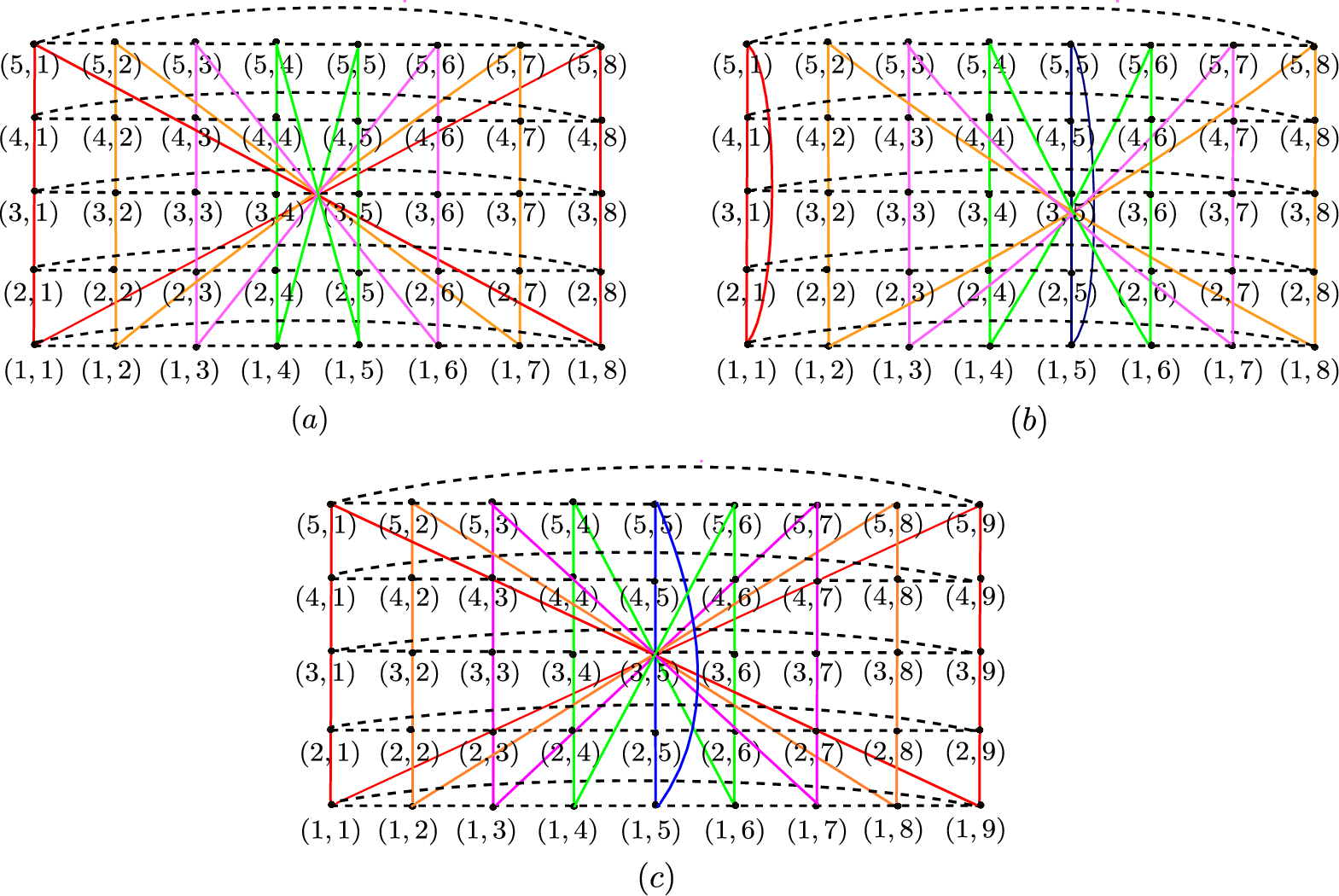}\\
\small Fig.3~~(a)~$C_5 \Box^{\varphi_{1}} C_8$, (b)~$C_5 \Box^{\varphi_{2}} C_8$, (c)~$C_5 \Box^{\varphi_{3}} C_9$, where $\varphi_{1}$ is a reflection without fixed points, $\varphi_{2}$ is a reflection with two fixed points, $\varphi_{3}$ is a reflection with one fixed point.
\end{figure}
\vspace{-1.2em}


\section{The case: $\varphi$ is a $d$-shift, $d \in \mathbb{Z}_{t}$}
The dispersability of $C_s \Box^{\varphi} C_t$ has been solved in $[4,8,9,10]$ if $\varphi$ is a trivial shift of $C_t$. By Lemma 2.4, in order to solve the dispersability of $C_s \Box^{\varphi} C_t$, it suffices to consider the case for $1 \leq d \leq \lfloor t/2 \rfloor $.
Based on the value of $gcd(t,d)$, we consider three cases to discuss the dispersability of $C_s \Box^{\varphi} C_t$.

\indent\textbf{Case 1: $gcd(t,d)=1$}\\
\indent If $gcd(t,d)=1$, we will show there exists an integer $k$ such that $C_s \Box^{\varphi} C_t$ is isomorphic to a circulant graph $C(\mathbb{Z}_{st},\{1,k\})$, whose dispersability has been completely solved in $[12]$.\\
\indent Here, we introduce a result about the Diophantine equation which is the main technique to calculate the position of an element in a special ordered set.\\
\noindent \textbf{Lemma 3.1.~}$^{[12]}$\emph{Given $a,b,c\in \mathbb{Z}^{+}$, if $gcd(a,b)=1$, then there exists a unique integer solution $(x_0,y_0)$ for the linear Diophantine equation $ax+by=c$ such that $0\leq x_0\leq b-1$. Furthermore, if $0\leq c \leq b-1$ and let $p(i)$ be the ordinal number which represents the position of the element $i$ in the ordered set $V_1^{\prime}=\{1,1+a,1+2a,\cdots,1+(b-1)a\}\pmod{b}$, then $p(1+c)$ is $1+x_0$.}\\
\noindent \textbf{Lemma 3.2.~}\emph{Let $G=C_s \Box^{\varphi} C_t$, where $\varphi$ is a $d$-shift$(1 \leq d \leq \lfloor t/2 \rfloor)$. If $gcd(t,d)=1$, let $(x_0,y_0)$ be the unique solution of the Diophantine equation $(t-d)x-ty=1$ such that $0\leq x_0\leq t-1$. Then $G\cong C(\mathbb{Z}_{st},\{1,sx_0\})$.}\\
\noindent \textbf{Proof.~}If $gcd(t,d)=1$, by Lemma 2.6, the spanning subgraph $G^{\prime}$ is an $st$-cycle. Next, we provide an algorithm to show the Cartesian graph bundle $G$ is a circulant graph $C(\mathbb{Z}_{st},\{1,sx_0\})$.\\
\noindent \textbf{Step 1.~Relabel the vertices of the Cartesian graph bundle $G$.}\\
\indent Let $D$ be an oriented cycle of $G^{\prime}$ and the vertex $(1,1)$ dominates the vertex $(2,1)$, etc. Since $G^{\prime}$ is a Hamilton cycle of $G$, relabel the vertices of $G$ along the cycle $G^{\prime}$ with the strictly increasing sequence $\{1,2,\cdots,st\}$. Without loss of generality, relabel the vertex $(1,1)$ with $1$, the vertex $(2,1)$ with $2$, etc. (See Fig.4 for the relabeling of $C_5\Box^{\varphi} C_7$, where $\varphi$ is a 3-shift.)

\vspace{-0.5em}
\begin{figure}[htbp]
\centering
\includegraphics[height=3.9cm, width=0.72\textwidth]{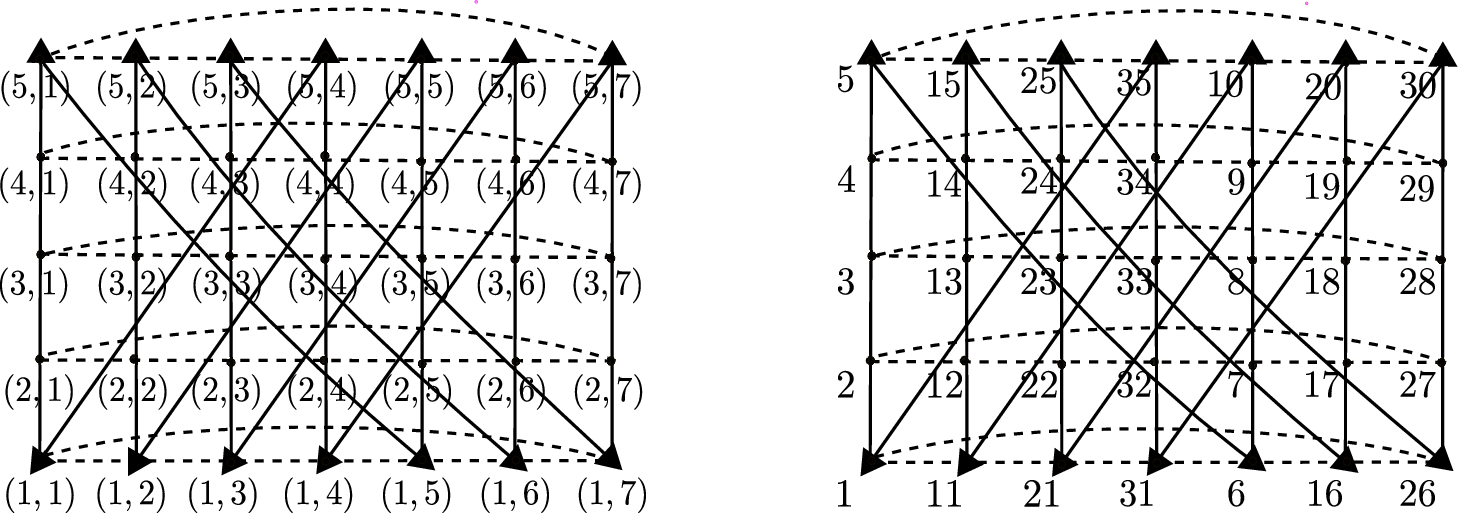}\\
\small{ Fig.4  ~Relabel the vertices of the Cartesian graph bundle $C_5\Box^{\varphi} C_7$ and the circulant graph $C(\mathbb{Z}_{35},\{1,10\})$~(right), where $\varphi$ is a 3-shift of $C_7$.}\\
\end{figure}
\vspace{-1.1em}

\noindent \textbf{Step 2.~Calculate the jump length of the circulant graph.}\\
\indent By the definition of circulant graph, observe that the Cartesian graph bundle $G$ is a circulant graph $C(\mathbb{Z}_{st},\{1,k\})$ for some $k$. Next we calculate the integer $k$. It is sufficient to calculate the integer $m$ which is the relabeling of the vertex $(1,2)$ and $k=m-1$.

\vspace{-0.5em}
\begin{figure}[htbp]
\centering
\includegraphics[height=3.9cm, width=0.7\textwidth]{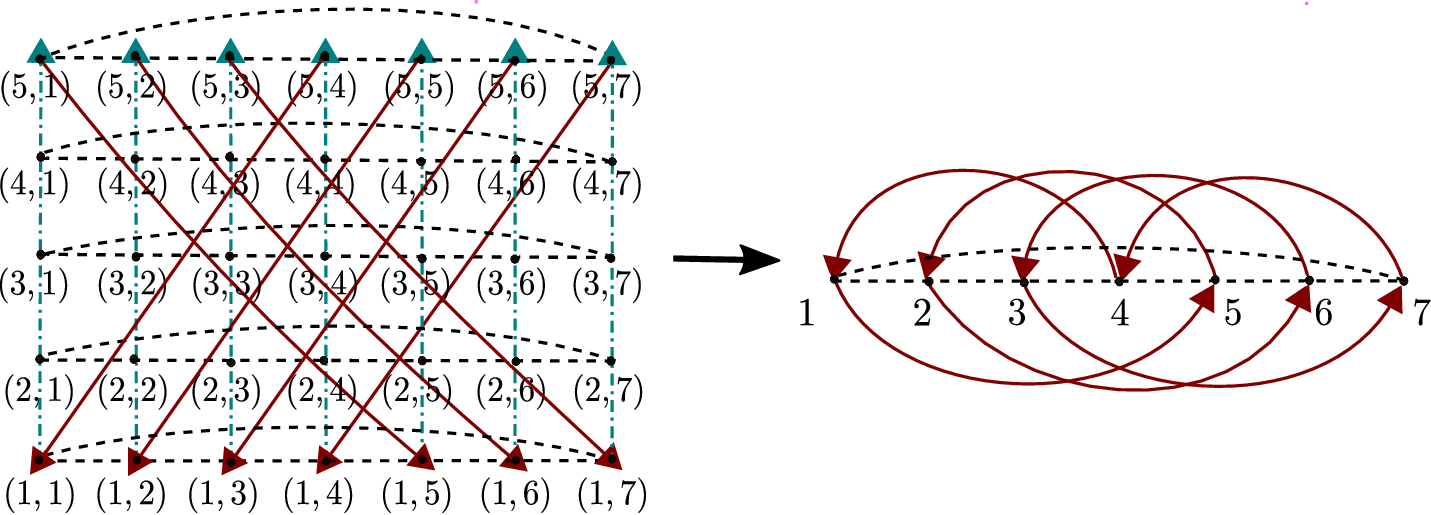}\\
\small{ Fig.5  ~The shrinking operation of the Cartesian graph bundle $C_5\Box^{\varphi} C_7$, where $\varphi$ is a 3-shift of $C_7$.}\\
\end{figure}
\vspace{-1.1em}

A graph $H^{\prime}$ can be obtained by the shrinking $B_j^{\prime}$ into a single vertex $j$, $j\in\mathbb{Z}_t$ (see Lemma 2.6). Let $D^{\prime}$ be an orientation of $H^{\prime}$ and keep the orientation of arcs which do not been contracted unchanged (see Fig.5 for $C_5\Box^{\varphi} C_7$). Therefore, $m=(p(2)-1)s+1$, where $p(2)$ is the ordinal number of the position of the element $2$ in the ordered set $\{1,1+(t-d),1+2(t-d),\cdots,1+(t-1)(t-d)\}(\text{mod}~t)$. Since $gcd(t,d)=1$, it is easy to see $gcd(t,t-d)=1$. By Lemma 3.1, that is to find the unique solution $(x_0,y_0)$ of the Diophantine equation $(t-d)x-ty=1$ such that $0 \leq x_0 \leq t-1$ and $m=x_0s+1$.
\hfill{$\square$}

\vspace{0.2em}
\textbf{Case 2: $gcd(t,d)$ is even}

\noindent \textbf{Lemma 3.3.~}\emph{Let $G=C_s \Box^{\varphi} C_t$, where $\varphi$ is a $d$-shift $(1 \leq d \leq \lfloor t/2 \rfloor)$. If $gcd(t,d)$ is even, then}
\vspace{-0.6em}
$$mbt(C_s \Box^{\varphi} C_t)=\left\{\begin{array}{l}
4,\ \ \ \ \ \ \ \ \ \ \ s~\text{is~even};\\
5,\ \ \ \ \ \ \ \ \ \ \ \ s~\text{is~odd}.\\
\end{array} \right.$$
\vspace{-1.8em}

\noindent \textbf{Proof.~}According to Lemma 2.1 - 2.3 and Fact 2.1, it is sufficient to show that 4 and 5 are the upper bounds of the matching book thickness of $G$ when $s$ is even and odd respectively. Put the vertices of $G$ clockwise along a circle in the order $V_1V_2^{-}V_3V_4^{-}\cdots V_{gcd(t,d)-1}V_{gcd(t,d)}^{-}$.~The edges of $G$ can be colored well in the following two steps (see Fig.6 for $C_6\Box^{\varphi}C_{10}$(left) and $C_5\Box^{\varphi}C_{12}$(right), where $\varphi$ is a 4-shift):

\vspace{-0.4em}
\begin{figure}[htbp]
\centering
\includegraphics[height=7.4cm, width=1\textwidth]{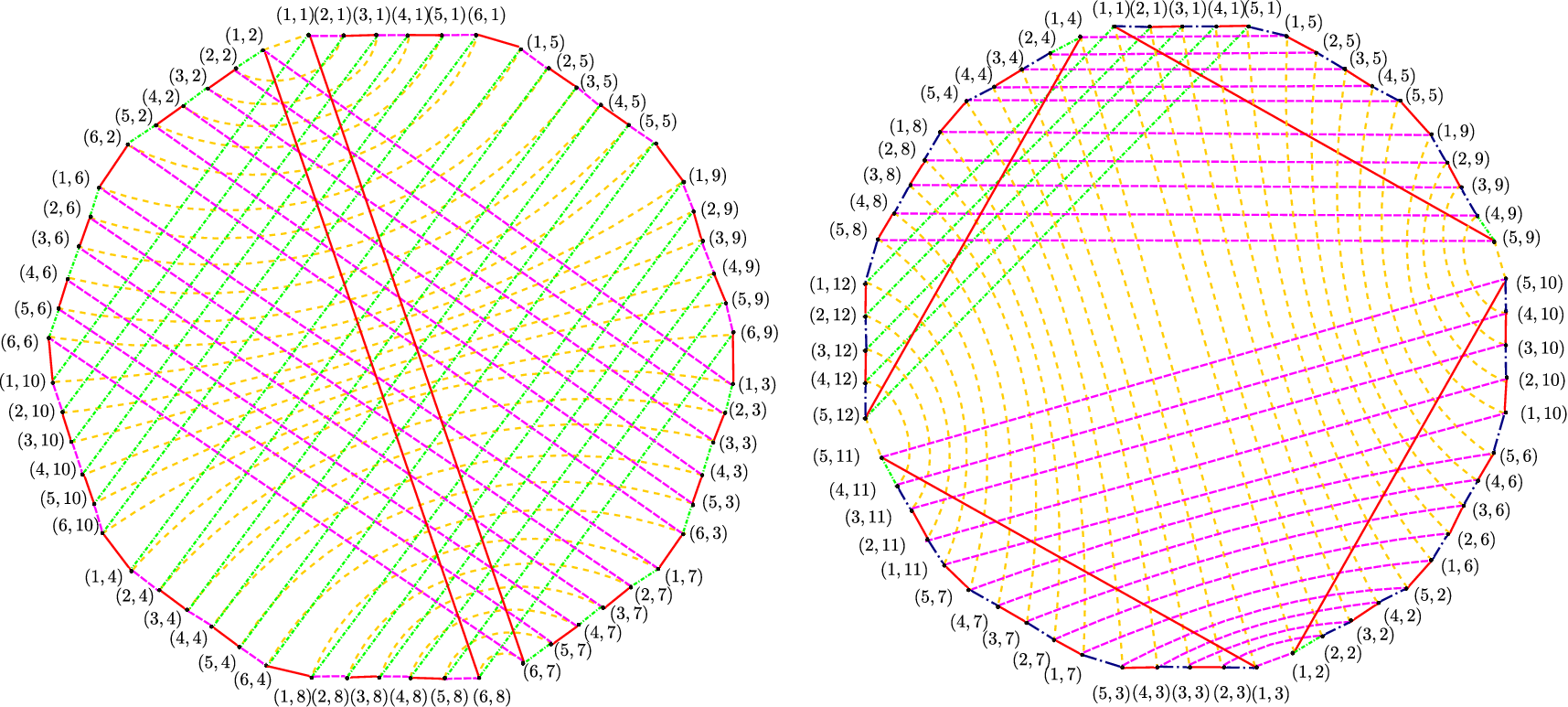}
\small Fig.6  ~A 4-page matching book embedding of $C_6\Box^{\varphi}C_{10}$(left) and a 5-page matching book embedding of $C_5\Box^{\varphi}C_{12}$(right), where $\varphi$ is a 4-shift.
\end{figure}
\vspace{-1.2em}

\noindent\textbf{Step 1:~}The coloring of the edges belonging to the cycles $F_1,F_2,\cdots,{F_s}$.\\
\indent Yellow: $\{((i,j),(i,j+1))~|~(i,j) \in V_m,m\in\{1,3,\cdots,gcd(t,d)-1\}\}$;\\
\indent Purple: $\{((i,j),(i,j+1))~|~(i,j) \in V_m,m\in\{2,4,\cdots,gcd(t,d)-2\}\}\cup\{((i,j),(i,j+1))~|~(i,j) \in V_{gcd(t,d)}, (i,j)\text{~is~before~the~element~}(1,t)\text{~in}~V_{gcd(t,d)}\}$;\\
\indent Green: $\{((i,j),(i,j+1))~|~(i,j) \in V_{gcd(t,d)}, (i,j)\text{~is~not~before~the~element~}(1,t)\text{~in}~V_{gcd(t,d)}\}$.\\
\noindent\textbf{Step 2:~}The coloring of the edges belonging to the cycles $H_1,H_2,\cdots,H_{gcd(t,d)}$.\\
\indent Color the edges of $\{((1,i),(s,i-d))~|~1\leq i \leq gcd(t,d)\}$ red. If $s$ is even, the cycles $H_1,H_2,\cdots,H_{gcd(t,d)}$ are even and red, green, purple can be used to color these even cycles well. Otherwise, use red, green and blue to color the cycles $H_1,H_2,\cdots,H_{gcd(t,d)}$ well.\\
\indent The result is established.
\hfill{$\square$}

\vspace{0.2em}
\textbf{Case 3: $gcd(t,d)\neq 1, gcd(t,d)$ is odd}\\
\noindent \textbf{Lemma 3.4.~}\emph{Let $G=C_s \Box^{\varphi} C_t$, where $\varphi$ is a $d$-shift $(1 \leq d \leq \lfloor t/2 \rfloor)$, $gcd(t,d)\neq1$ and $gcd(t,d)$ is odd. We have}

\vspace{-1.2em}
$$mbt(C_s \Box^{\varphi} C_t)=\left\{\begin{array}{l}
4,\ \ \ \ t~\text{is~even},~s,d~\text{are~both~odd};\\
5,\ \ \ \ \ \ \ \ \ \ \ \ \ \ \ \ \ \ \ \ \ \ \ \ \ \ \ \text{otherwise}.\\
\end{array} \right.$$
\vspace{-1.4em}

\noindent \textbf{Proof.~}According to Lemma 2.3, the graph $G$ is bipartite if $t$ is even and $s,d$ are both odd; otherwise, $G$ is nonbipartite. We consider the dispersability of $G$ based on whether $G$ is bipartite.

\indent\textbf{Case 1. $G$ is bipartite}\\
\indent By Lemma 2.2 and Fact 2.1, it is sufficient to show $mbt(G)\leq 4$. Put the vertices of $G$ clockwise along a circle in the order $B_{1}B_{t}^{-}B_{3}B_{t-2}^{-}\cdots B_{t-1}^{}B_{2}^{-}$.~The edges of $G$ can be colored well in the following two steps (see Fig.7 for $C_5\Box^{\varphi}C_{20}$, where $\varphi$ is a 5-shift):

\vspace{-0.5em}
\begin{figure}[htbp]
\centering
\includegraphics[height=10cm, width=0.62\textwidth]{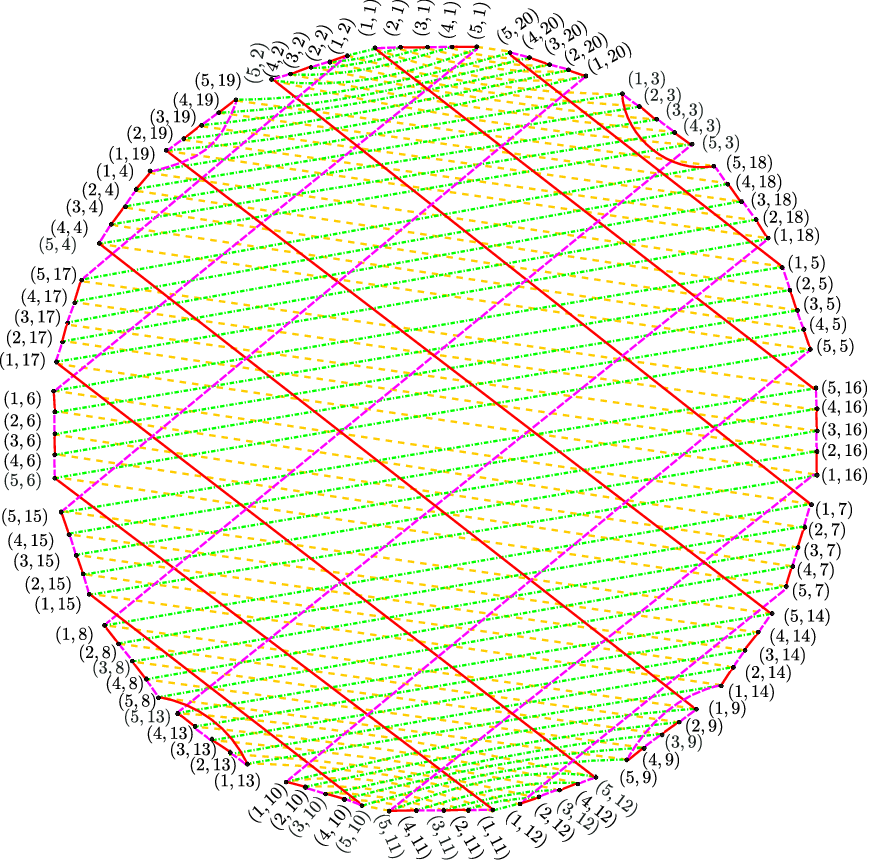}
\small\centerline{Fig.7  ~A 4-page matching book embedding of $C_5\Box^{\varphi}C_{20}$, where $\varphi$ is a 5-shift of $C_{20}$.}%
\end{figure}
\vspace{-1.3em}

\noindent\textbf{Step 1:~}The coloring of the edges belonging to the cycles $F_1,F_2,\cdots,{F_s}$.\\
\indent Yellow: $\{((i,j),(i,j+1))~|~1\leq i\leq s, j \in\{2,4,\cdots,t\}\}$;\\
\indent Green:~$\{((i,j),(i,j+1))~|~1\leq i\leq s, j\in\{1,3,\cdots,t-1\}\}$.\\
\noindent\textbf{Step 2:~}The coloring of the edges belonging to the cycles $H_1,H_2,\cdots,H_{gcd(t,d)}$.\\
\indent Color the edges of $\{((1,i),(s,i-d))~|~i\in\{1,3,\cdots,t-1\}\}$ red,~the edges of $\{((1,i),(s,i-d))~|~\{2,4,\cdots,t\}\}$ purple.~The remaining even edges of the path induced by $B_i^{\prime}$ can be colored well with red and purple alternately, $i=1,2,\cdots,t$.\\
\indent\textbf{Case 2. $G$ is nonbipartite}\\
\indent According to Lemma 2.1 - 2.2 and Fact 2.1, it is sufficient to show $mbt(G)\leq 5$. Let ordered set $X$ be a rearrangement of $V_{1}\cup V_{2}$ and we denote the $jth$ element of an ordered vertex set $V_i$ by $V_{i_j}$. If $|V_1|$ is even, let $X=\{V_{1_1},V_{2_1},V_{2_2},V_{1_2},V_{1_3},V_{2_3},V_{2_4},V_{1_4},\cdots,V_{1_{|V_1|-1}},V_{2_{|V_1|-1}},V_{2_{|V_1|}},V_{1_{|V_1|}}\}$; otherwise, let $X=\{V_{1_1},V_{2_1},V_{2_2},V_{1_2},V_{1_3},V_{2_3},V_{2_4},V_{1_4},\cdots,V_{1_{|V_1|-2}},V_{2_{|V_1|-2}},V_{2_{|V_1|-1}},V_{1_{|V_1|-1}},V_{1_{|V_1|}}, V_{2_{|V_1|}}\}$. Put the vertices of $G$ clockwise along a circle in the order $XV_{3}^{-}V_{4}V_5^{-}\cdots V_{gcd(t,d)-1}V_{gcd(t,d)}^{-}$. The edges of $G$ can be colored well in the following two steps (see Fig.8 for $C_4\Box^{\varphi}C_{15}$(left) and $C_5\Box^{\varphi}C_{15}$(right), where $\varphi$ is a 5-shift):\\
\noindent\textbf{Step 1:~}The coloring of the edges belonging to the cycles $F_1,F_2,\cdots,{F_s}$.\\
\indent\textbf{Subcase 2.1 $H_1$ is an even cycle}\\
\indent Green: $\{((i,j),(i,j+1))~|~(i,j) \in V_k, k\in\{2,4,\cdots,gcd(t,d)-1\}\}$;\\
\indent Yellow: $\{((i,j),(i,j+1))~|~(i,j) \in V_k, k\in\{3,5,\cdots,gcd(t,d)-2\}\}\cup\{((i,j),(i,j+1))~|~(i,j) \in V_{gcd(t,d)}, j\neq t \}\cup\{((i,j),(i,j+1))~|~(i,j)\in V_1,j= 1\}$;\\
\indent Purple: $\{((i,j),(i,j+1))~|~(i,j) \in V_{gcd(t,d)}, j= t \}\cup\{((i,j),(i,j+1))~|~(i,j)\in V_1,j\neq 1\}$.
\indent\textbf{Subcase 2.2 $H_1$ is an odd cycle}\\
\indent Yellow: $\{((i,j),(i,j+1))~|~(i,j) \in V_k, k\in\{3,5,\cdots,gcd(t,d)-2\}\}\cup\{((i,j),(i,j+1))~|~(i,j) \in V_{gcd(t,d)}, j\neq t \}\cup\{((i,1),(i,2))~|~3\leq i\leq s\}$;\\
\indent Green: $\{((i,j),(i,j+1))~|~(i,j) \in V_k, k\in\{2,4,\cdots,gcd(t,d)-1\}, (i,j)\neq (s,2-d)\}\cup\{((s,1-d),(s,2-d))\}$;\\
\indent Blue: $\{((s,2-d),(s,3-d))\}\cup\{((2,1),(2,2))\}$;\\
\indent Purple: $\{((i,t),(i,1))~|~2\leq i \leq s\}\cup\{((1,1),(1,2))\}\cup\{((i,j),(i,j+1))~|~(i,j)\in V_1,j\neq 1,(i,j)\notin \{(s-1,1-d),(s,1-d)\}\}$;\\
\indent Red: $\{((1,1),(1,t))\}\cup\{((s-1,1-d),(s-1,2-d))\}$.

\vspace{-0.9em}
\begin{figure}[htbp]
\centering
\includegraphics[height=5.9cm, width=1\textwidth]{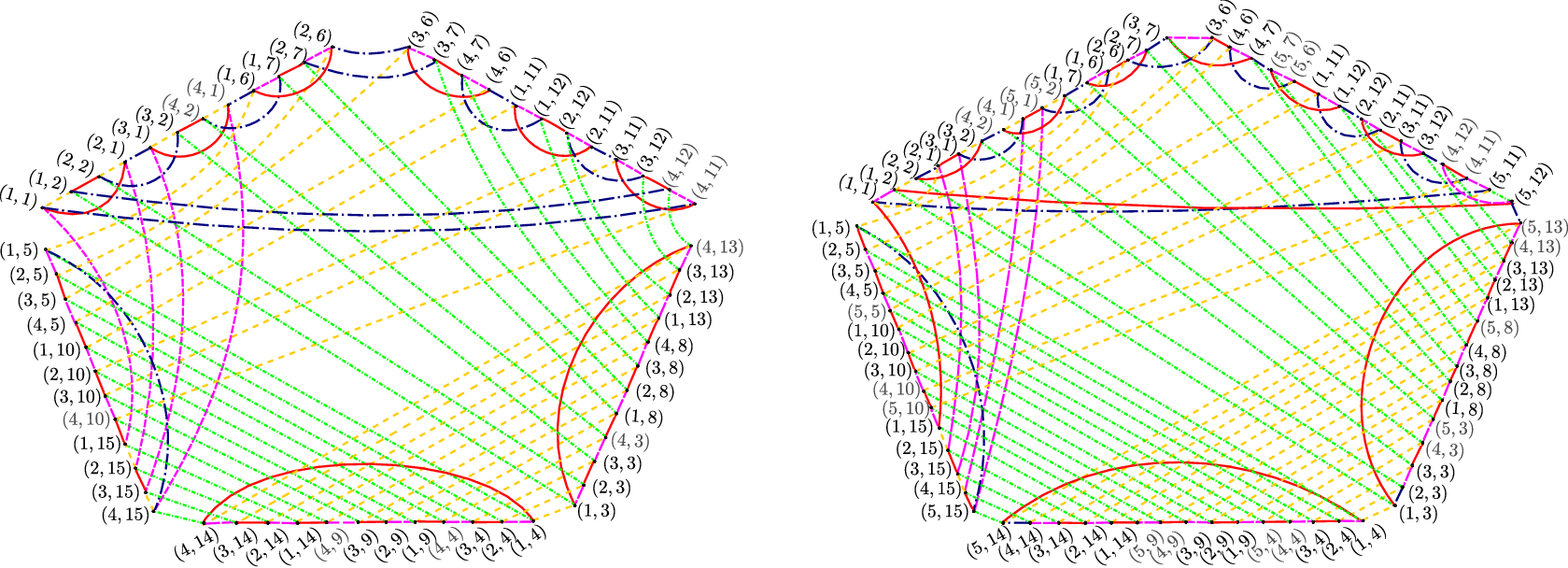}
\small Fig.8  ~A 5-page matching book embedding of $C_4\Box^{\varphi}C_{15}$(left) and a 5-page matching book embedding of $C_5\Box^{\varphi}C_{15}$(right), where $\varphi$ is a 5-shift of $C_{15}$.
\end{figure}
\vspace{-1.2em}

\noindent\textbf{Step 2:~}The coloring of the edges belonging to the cycles $H_1,H_2,\cdots,H_{gcd(t,d)}$.\\
\indent Color the cycles $H_3,H_4,\cdots,H_{gcd(t,d)-1}$ well with red, purple and blue.~As for the edges of the cycle $H_1$ and $H_2$, if the cycle $H_1$ is an even cycle, color the edges of $H_1$ and $H_2$ with red and blue alternately, where $((1,1),(2,1))$, $((1,2),(2,2))$ are red. Otherwise, color the edges $((1,1),(2,1))$, $((1,2),(2,2))$ yellow, the edges $((s-1,2-d),(s,2-d))$, $((s-1,1-d),(s,1-d))$ purple, the edge $((1,1),(s,1-d))$ blue, the edge $((1,2),(s,2-d))$ red, and the remaining even edges with red and blue alternately, where $((2,1),(3,1))$, $((2,2),(3,2))$ are red. For the cycle $H_{gcd(t,d)}$, color the edge $((1,gcd(t,d)),(s,t))$ blue, the edges of the path induced by $\{(1,t),(2,t),\cdots,(s,t)\}$ well with yellow and red alternately, the remaining edges well with purple and red alternately.\\
\indent Therefore, the result is established.
\hfill{$\square$}

\vspace{0.2em}
In conclusion, the following conclusion holds.

\noindent\textbf{Theorem 3.1.~}\emph{Let $G=C_s \Box^{\varphi} C_t$, where $s,t\geq 3$ and $\varphi$ is a $d$-shift, $d\in \mathbb{Z}_t$. We have}
\vspace{-0.5em}
$$mbt(G)=\left\{\begin{array}{l}
4,\ \ \ \ \ \ \ \ \ \ G~$is~bipartite$;\\
5,\ \ \ \ \ \ \ \ \ \ \ \ \ \ \ \ $otherwise$.\\
\end{array} \right.$$


\section{The case: $\varphi$ is a reflection}
According to the parity of $s$, we discuss the dispersability of $C_s \Box^{\varphi} C_t(s,t\geq 3)$ in the following two cases, where $\varphi$ is a reflection of $C_t$.

\vspace{0.2em}
\textbf{Case 1: $s$ is odd}\\
\noindent\textbf{Lemma 4.1.~}\emph{Let $G=C_s \Box^{\varphi} C_t$, where $s$ is odd and $\varphi$ is a reflection. We have}
\vspace{-0.6em}
$$mbt(C_s \Box^{\varphi} C_t)=\left\{\begin{array}{l}
4,\ \ \ \ \ \ \ \varphi~$has~no~fixed~points$;\\
5,\ \ \ \ \ \ \ \ \ \ \ \ \ \ \ \ \ \ \ \ \ \ \ $otherwise$.\\
\end{array} \right.$$
\vspace{-1.9em}

\noindent \textbf{Proof.}
According to Lemma 2.3, the graph $G$ is bipartite if $\varphi$ has no fixed points; otherwise, $G$ is nonbipartite. By Lemma 2.1 - 2.2 and Fact 2.1, it is sufficient to show that 4 and 5 are the upper bounds of the matching book thickness of $G$ when $G$ is bipartite and nonbipartite respectively. Put the vertices of $G$ clockwise along a circle in the order $A_1^{-}A_2A_3^{-}A_4\cdots A_{s-1}A_{s}^{-}$.~The edges of $G$ can be colored well in the following two steps:

\vspace{-0.9em}
\begin{figure}[htbp]
\centering
\includegraphics[height=5.8cm, width=1\textwidth]{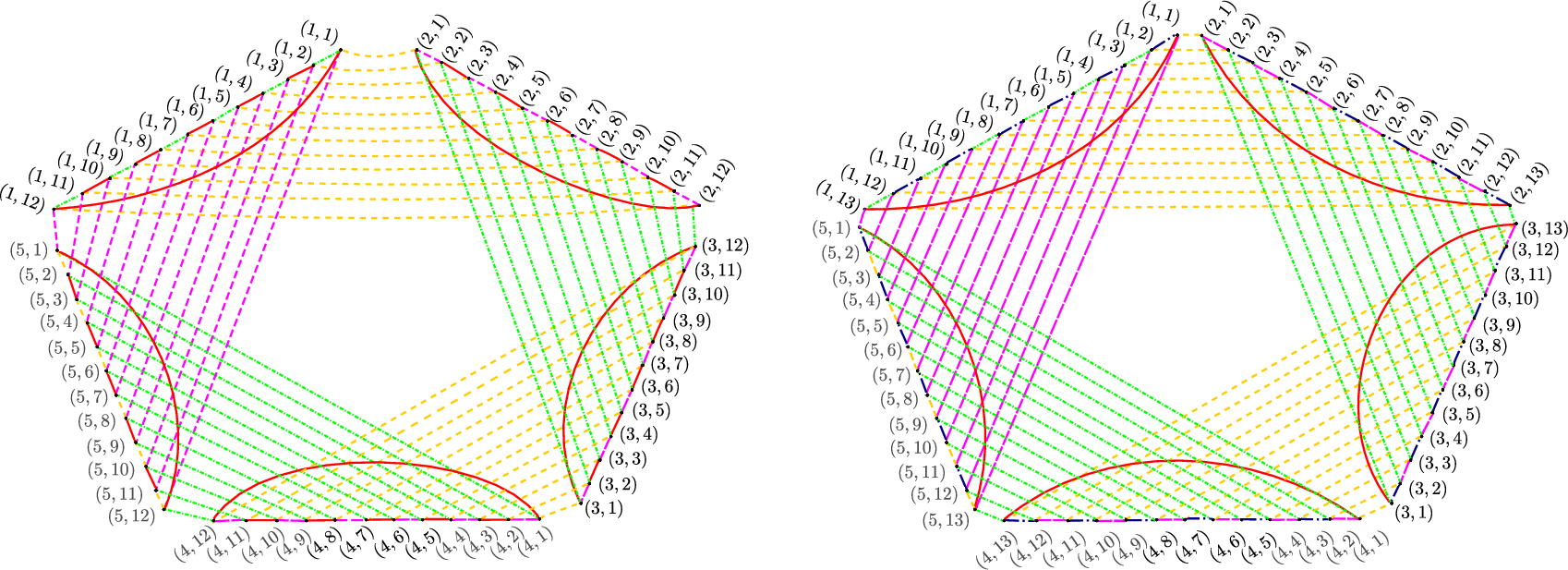}\\
\small Fig.9 ~A 4-page matching book embedding of $C_{5} \Box^{\varphi_{1}} C_{12}$(left), where $\varphi_{1}$ is a reflection without fixed points and $C_5 \Box^{\varphi_{2}}C_{13}$(right), where $\varphi_{2}$ is a reflection with one fixed point.
\end{figure}
\vspace{-1.3em}

\noindent\textbf{Step 1:~}The coloring of the edges not belonging to the cycles $F_1,F_2,\cdots,{F_s}$.\\
\indent Color the edges of $\{((i,j),(i+1,j))~|~i \in \{1,3,\cdots,s-2\}, j \in \{1,2,\cdots,t\} \}$ yellow, the edges of $\{((i,j),(i+1,j))~|~i \in \{2,4,\cdots,s-1\}, j \in \{1,2,\cdots,t\} \}$ green. If $\varphi$ has 0 or 1 fixed point, color the edges of $\{((s,i),(1,t+1-i))~|~1 \leq i \leq t \}$ purple (see Fig.9 for the matching book embedding of $C_{5} \Box^{\varphi_1} C_{12}$(left), where $\varphi_1$ has no fixed points and $C_5 \Box^{\varphi_2} C_{13}$(right), where $\varphi_2$ has one fixed point). Otherwise, color the edges of $\{((1,i),(s,t+2-i))~|~2 \leq i \leq t \}$ purple, the edge $((1,1),(s,1))$ blue (see Fig.10 for $C_s \Box^{\varphi_3}C_t$(left) and a 5-page matching book embedding of $C_{5} \Box^{\varphi_3} C_{12}$(right), where $\varphi_3$ has two fixed points).\\
\noindent\textbf{Step 2:~}The coloring of the edges belonging to the cycles $F_1,F_2,\cdots,{F_s}$.\\
\indent It is easy to use yellow, green, purple, red to color the cycles $F_1,F_2,\cdots,{F_s}$ well if $t$ is even; use yellow, green, purple, red, blue to color the cycles $F_1,F_2,\cdots,{F_s}$ well if $t$ is odd.

\vspace{-0.3em}
\begin{figure}[htbp]
\centering
\includegraphics[height=6.2cm,width=1\textwidth]{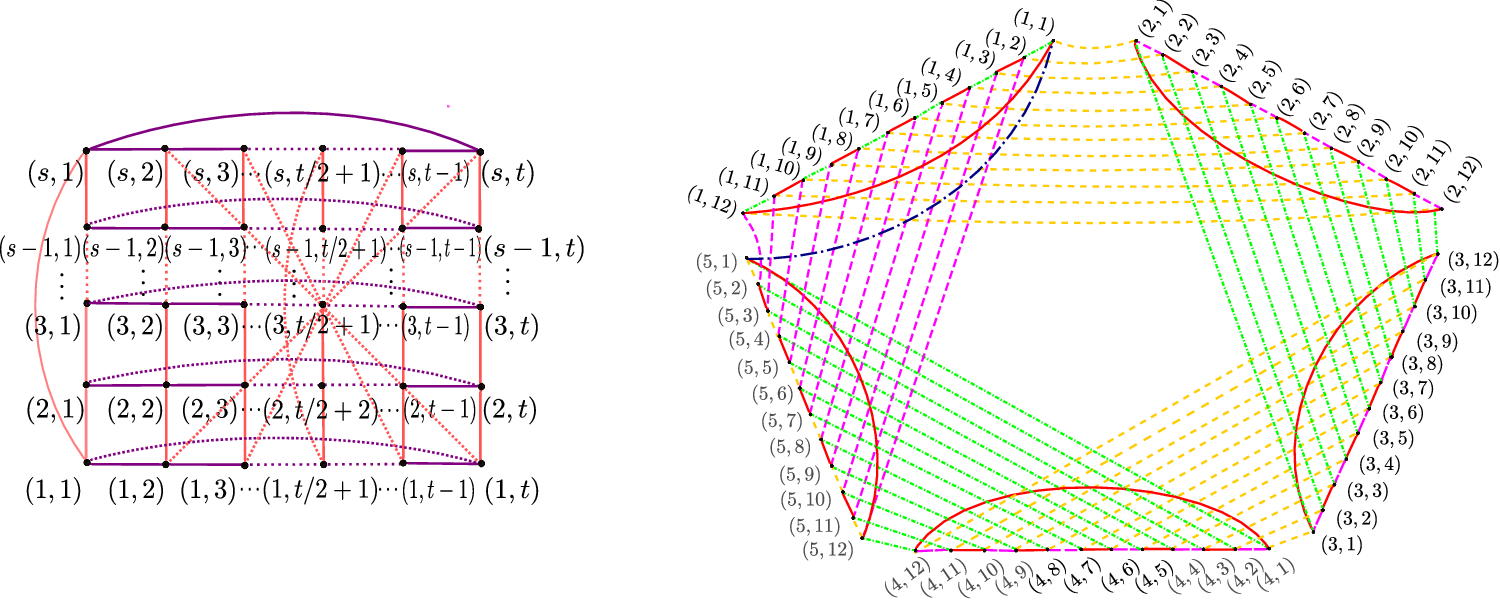}
\small Fig.10~~$C_s \Box^{\varphi_3}C_t$(left) and a 5-page matching book embedding of $C_{5}\Box^{\varphi_3} C_{12}$(right), where $\varphi_3$ has two fixed points.
\end{figure}
\vspace{-1.2em}

The result is established.
\hfill{$\square$}

\vspace{0.2em}
\textbf{Case 2: $s$ is even}\\
\noindent \textbf{Lemma 4.2.~}\emph{Let $G=C_s \Box^{\varphi} C_t$, where $s$ is even and $\varphi$ is a reflection. We have }\\
\vspace{-0.8em}
$$mbt(C_s \Box^{\varphi} C_t)=\left\{\begin{array}{l}
4,\ \ \ \ \ \ \ \ \varphi~$has~two~fixed~points$;\\
5,\ \ \ \ \ \ \ \ \ \ \ \ \ \ \ \ \ \ \ \ \ \ \ \ \ $otherwise$.\\
\end{array} \right.$$

\noindent \textbf{Proof.~}
\textbf{Case 1: $\varphi$ has two fixed points.}

According to Lemma 2.2 and Fact 2.1, it suffices to show that $mbt(G)\leq 4$. Put the vertices of $G$ clockwise along a circle in the order $B_{1}B_{2}^{-}B_{3}B_{4}^{-}\cdots B_{t-1}B_{t}^{-}$.~The edges of $G$ can be colored well in the following two steps (see Fig.11(left) for $C_6 \Box^{\varphi_1}C_{8}$, where $\varphi_1$ has two fixed points):

\noindent\textbf{Step 1:~}The coloring of the edges belonging to the cycles $F_1,F_2,\cdots,{F_s}$.\\
\indent Yellow: $\{((i,j),(i,j+1))~|~1\leq i \leq s,j \in \{1,3,\cdots,t-1\}\}$;\\
\indent Green: $\{((i,j),(i,j+1))~|~1\leq i \leq s,j \in \{2,4,\cdots,t\}\}$.\\
\noindent\textbf{Step 2:~}The coloring of the edges not belonging to the cycles $F_1,F_2,\cdots,{F_s}$.\\
\indent By Lemma 2.7, color all remaining even cycles well with blue and red alternately.

\textbf{Case 2: $\varphi$ has one fixed point.}\\
\indent According to Lemma 2.1 - 2.3 and Fact 2.1, it is sufficient to show that $mbt(G)\leq 5$. Let $U$ be a rearrangement of $B_1\cup B_t$ and $U=\{(s,1),(s,t),(s-1,t),(s-1,1),(s-2,1),(s-2,t),(s-3,t),(s-3,1),\cdots,(4,1),(4,t),(3,t),(3,1),(2,1),(2,t),(1,1),(1,t)\}$. Put the vertices of $G$ clockwise along a circle in the order $U^{-}B_2^{-}B_3B_4^{-}B_5^{}\cdots B_{t-2}B_{t-1}^{-}$. The edges of $G$ can be colored well in the following two steps (see Fig.11(right) for $C_6 \Box^{\varphi_2}C_{9}$, where $\varphi_2$ has one fixed point):\\
\noindent\textbf{Step 1:~}The coloring of the edges belonging to the cycles $F_1,F_2,\cdots,{F_s}$.\\
\indent Yellow: $\{((i,j),(i,j+1))~|~1\leq i\leq s, j\in\{1,3,\cdots, t-2\},(i,j)\neq (2,1)\}\cup\{((2,1),(2,t))\}$;\\
\indent Green: $\{((i,j),(i,j+1))~|~1\leq i\leq s, j\in\{2,4,\cdots,t-1\},(i,j)\neq (1,t-1)\}\cup\{((1,1),(1,t))\}$;\\
\indent Blue: $\{((1,t-1),(1,t))\}$;\\
\indent Red: $\{((2,1),(2,2))\}\cup\{((i,1),(i,t))~|~3\leq i\leq s\}$.\\
\noindent\textbf{Step 2:~}The coloring of the edges not belonging to the cycles $F_1,F_2,\cdots,{F_s}$.\\
\indent Color the even cycles induced by $B_i\cup B_{t+1-i}(3\leq i\leq (t-1)/2),~B_{(t+1)/2}$ with blue and purple alternately.~For $i\in\{1,2\}$, color the edges $((1,t+1-i),(2,t+1-i))$ red,~the remaining edges induced by $B_i\cup B_{t+1-i}$ with blue and purple alternately, where $((1,t+1-i),(s,i))$ are purple.

\vspace{-0.3em}
\begin{figure}[htbp]
\centering
\includegraphics[height=7.2cm,width=1\textwidth]{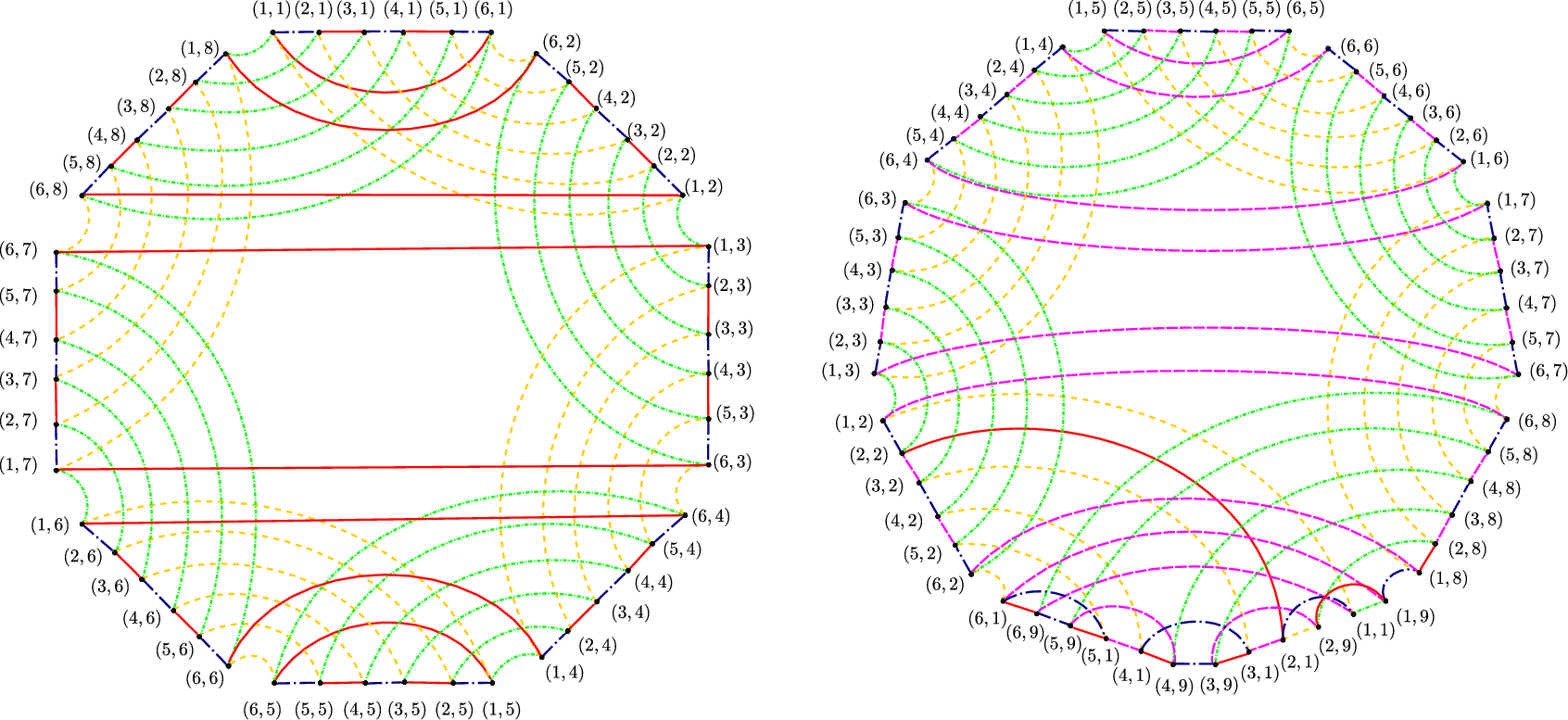}\\
\small Fig.11~The matching book embedding of $C_6 \Box^{\varphi_1}C_{8}$(left), where $\varphi_1$ has two fixed points and $C_6 \Box^{\varphi_2}C_{9}$(right), where $\varphi_2$ has one fixed point.
\end{figure}
\vspace{-1.3em}

\vspace{0.2em}
\textbf{Case 3: $\varphi$ has no fixed points.}\\
\indent By Lemma 2.1 - 2.3 and Fact 2.1, it suffices to show that $mbt(G)\leq 5$. Let $V_i$ and $W_i$ be two ordered subsets of $B_i$, specifically, $V_i=\{(1,i),(s,i)\}$, $W_i=B_i\setminus V_i$, $i=1,2,\cdots,t$. Assume ordered sets $M_j=V_j\cup V_{t+1-j}$, $N_j=W_j\cup W_{t+1-j}^{-}$, $j=1,2,\cdots,t/2$. If $t/2$ is odd, put the vertices of $G$ clockwise along a circle in the order $M_1M_2^{-}M_3M_4^{-}\cdots M_{t/2-1}^{-}M_{t/2}N_{t/2}^{-}N_{t/2-1}N^{-}_{t/2-2}N_{t/2-3}\cdots N^{-}_3 N_2\\N_1^{-}$ (see Fig.12a for $C_8 \Box^{\varphi}C_{10}$); otherwise, along a circle in the order $M_1M_2^{-}M_3M_4^{-}\cdots M_{t/2-1}M_{t/2}^{-}\\N_{t/2}N^{-}_{t/2-1}N_{t/2-2}N^{-}_{t/2-3}\cdots N_2N_1^{-}$ (see Fig.12b for $C_8 \Box^{\varphi}C_8$). The edges of $G$ can be colored well in the following two steps:

\noindent\textbf{Step 1:~}The coloring of the edges belonging to the cycles $F_1,F_2,\cdots,{F_s}$.

Green: $\{((i,j),(i,j+1))~|~(i,j)\in B_2\cup B_4 \cup \cdots\cup B_{t-2} \cup B_t \setminus \{V_{t/2}\cup V_t\}\}$;

Yellow: $\{((i,j),(i,j+1))~|~(i,j)\in B_1\cup B_3 \cup \cdots\cup B_{t-3} \cup B_{t-1} \setminus V_{t/2}\}$;

Blue: $\{((1,i),(1,i+1))~|~i\in\{t/2,t\}\}$;

Red: $\{((s,i),(s,i+1))~|~i\in\{t/2,t\}\}$.

\noindent\textbf{Step 2:~}The coloring of the edges not belonging to the cycles $F_1,F_2,\cdots,{F_s}$.\\
\indent Color the edges of $\{((s-1,i),(s,i))~|~i\in\{1,t/2\}\}\cup\{((1,j),(2,j))~|~j\in\{t/2+1,t\}\}$ purple, the edges of $\{((s-1,i),(s,i))~|~i\in\{t/2+1,t\}\}$ blue, the edges of $\{((1,i),(2,i))~|~i\in\{1,t/2\}\}$ red, the edges of $\{((i,j),(s+1-i,t+1-j))~|~(i,j)\in V_1\}$ green, the edges of $\{((i,j),(s+1-i,t+1-j))~|~(i,j)\in V_2\cup V_3\cup\cdots\cup V_{t/2-1}\}$ red. If $t/2$ is odd, color the edges of $\{((i,j),(s+1-i,t+1-j))~|~(i,j)\in V_{t/2}\}$ yellow; otherwise, color them green. For the remaining edges of the paths induced by $B_i(1\leq i\leq t)$, use red, purple and blue to color them well.

This completes the proof.
\hfill{$\square$}

\vspace{-0.5em}
\begin{figure}[htbp]
\centering
\includegraphics[height=12cm,width=0.72\textwidth]{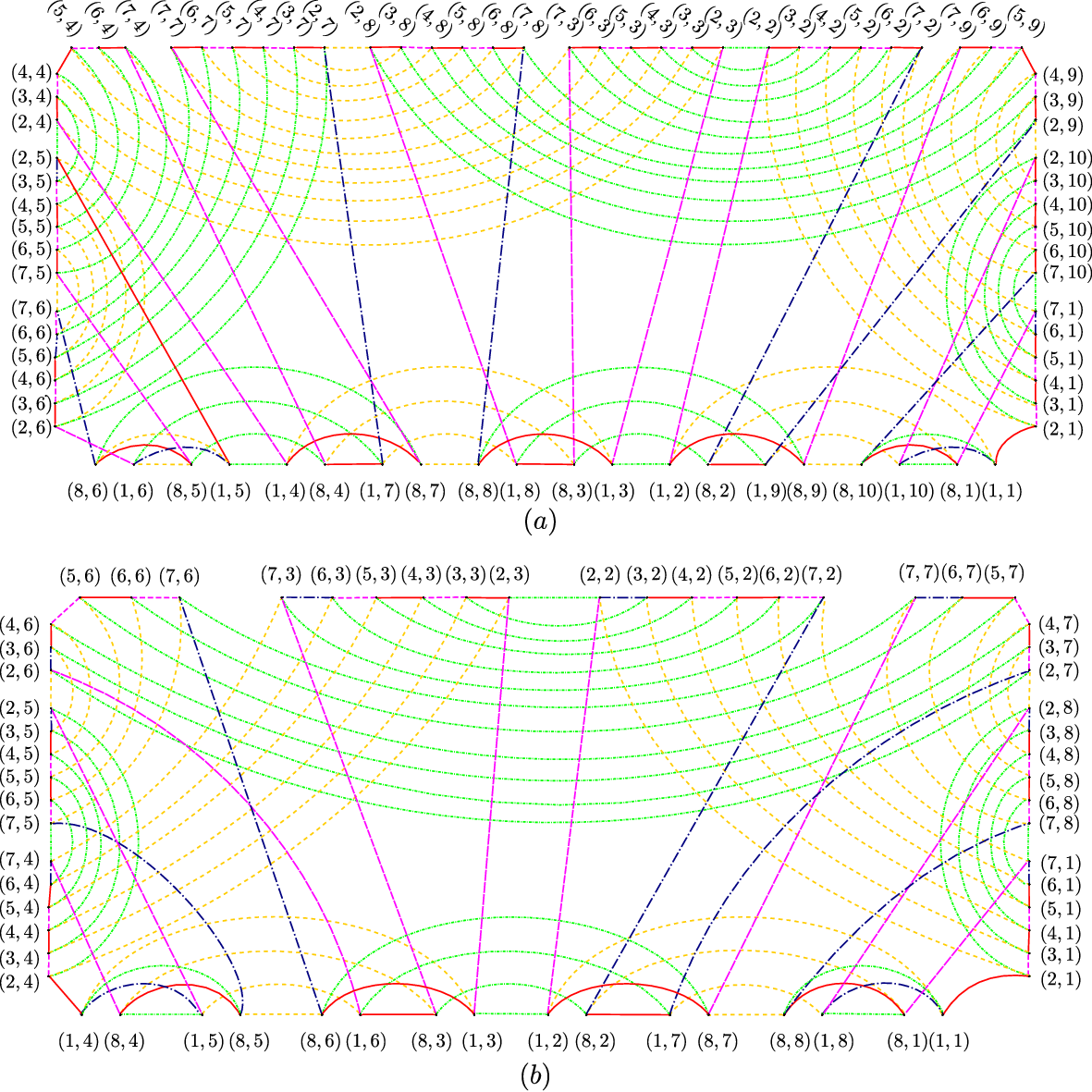}\\
\small Fig.12~A 5-page matching book embedding of $C_8 \Box^{\varphi}C_{10}$($a$) and $C_8 \Box^{\varphi}C_{8}$($b$), where $\varphi$ has no fixed points.
\end{figure}
\vspace{-1em}

By Lemma 4.1 and Lemma 4.2, the following conclusion holds.

\noindent\textbf{Theorem 4.1.~}\emph{Let $G=C_s \Box^{\varphi} C_t$, where $s,t\geq 3$ and $\varphi$ is a reflection of $C_t$. We have}
\vspace{-0.5em}
$$mbt(G)=\left\{\begin{array}{l}
4,\ \ \ \ \ \ \ \ \ \ G~$is~bipartite$;\\
5,\ \ \ \ \ \ \ \ \ \ \ \ \ \ \ \ $otherwise$.\\
\end{array} \right.$$

\section{Conclusion}
\indent

In this work, the dispersability of the Cartesian graph bundle of two cycles is completely solved , which implies the matching book thickness of the circulant graph $C(\mathbb{Z}_{n},\{k_1,k_2\})$ with two jump lengths is completely solved. Specifically, the graph $C(\mathbb{Z}_{n},\{k_1,k_2\})$ are dispersable if they are bipartite; otherwise, they are nearly dispersable. In $[15]$, Alam etc. ask whether vertex-transitive regular bipartite graphs are dispersable. As is well known, connected circulant graphs are vertex-transitive graphs. Our work partially supports the positive side.\\

\section*{Acknowledgment}

\indent

This work was partially funded by Science and Technology Project of Hebei Education Department, China (No. ZD2020130) and the Natural Science Foundation of Hebei Province, China (No. A2021202013).


%
%
%
%
%


\end{document}